\title{An application of type sequences to the blowing-up.}

\documentclass{article}   

 \author{Anna Oneto and Elsa Zatini\\ (Beitrage zur  Algebra und Geometrie, Vol 46, n.2, pp 471-489 (2005))}
\date{}

\expandafter\chardef\csname pre amssym.def  at\endcsname=\the\catcode`\@

\catcode`\@=11

\def\undefine#1{\let#1\undefined}
\def\newsymbol#1#2#3#4#5{\let\next@\relax
 \ifnum#2=\@ne\let\next@\msafam@\else
 \ifnum#2=\tw@\let\next@\msbfam@\fi\fi
 \mathchardef#1="#3\next@#4#5}
\def\mathhexbox@#1#2#3{\relax
 \ifmmode\mathpalette{}{\m@th\mathchar"#1#2#3}%
 \else\leavevmode\hbox{$\m@th\mathchar"#1#2#3$}\fi}
\def\hexnumber@#1{\ifcase#1 0\or 1\or 2\or 3\or 4\or 5\or 6\or 7\or 8\or
 9\or A\or B\or C\or D\or E\or F\fi}

\font\tenmsa=msam10
\font\sevenmsa=msam7
\font\fivemsa=msam5
\newfam\msafam
\textfont\msafam=\tenmsa
\scriptfont\msafam=\sevenmsa
\scriptscriptfont\msafam=\fivemsa
\edef\msafam@{\hexnumber@\msafam}

\font\tenmsb=msbm10
\font\sevenmsb=msbm7
\font\fivemsb=msbm5
\newfam\msbfam
\textfont\msbfam=\tenmsb
\scriptfont\msbfam=\sevenmsb
\scriptscriptfont\msbfam=\fivemsb
\edef\msbfam@{\hexnumber@\msbfam}

   \font\tengothic=eufm10
   \font\sevengothic=eufm7
   \newfam\gothicfam
   \textfont\gothicfam=\tengothic
   \scriptfont\gothicfam=\sevengothic
   \def\goth#1{{\fam\gothicfam #1}}
   \font\tenmsb=msbm10
   \font\sevenmsb=msbm7
   \newfam\msbfam
   \textfont\msbfam=\tenmsb
   \scriptfont\msbfam=\sevenmsb
   
\newtheorem{prop}{Proposition}[section]
\newtheorem{rem}[prop]{Remark}
\newtheorem{thm}[prop]{Theorem}
\newtheorem{coro}[prop]{Corollary}
\newtheorem{defin}[prop]{Definition}
\newtheorem{lemma}[prop]{Lemma}

\newtheorem{example}[prop]{Example}
\newtheorem{(*)}[prop]{}

\newcommand{\Ib}{\mbox{$\overline{I}$}}
\newcommand{\It}{\mbox{$\widetilde{I}$}}
\newcommand{\cx}{\mbox{${\rm C\hspace{-1.7mm}\sf l\hspace{1mm}}$}}

\newcommand{\integ}{\mbox{${\sf Z\hspace{-1.3 mm}Z}$}}

\newcommand{\nat}{\mbox{${\rm I\hspace{-.6 mm}N}$}}
\newcommand{\R}{\mbox{$\overline{R}$}}

\newcommand{\m}{\mbox{${\goth m}$}}

\newcommand{\w}{\mbox{$\omega$}}

\newcommand{\LL}{\mbox{$\Lambda$}}

\font\tengothic=eufm10
\font\sevengothic=eufm7
\newfam\gothicfam
\textfont\gothicfam=\tengothic
\scriptfont\gothicfam=\sevengothic
\def\goth#1{{\fam\gothicfam #1}}

\begin{document}
 
\maketitle

\noindent {\bf Abstract} \\
{\it  Let I be an $\m$-primary ideal of a
one-dimensional, analytically irreducible and 
residually rational local Noetherian
 domain  R. Given the   blowing-up  of R 
   along   $ I,$  we establish connections between  the type-sequence of  R and  classical invariants  like multiplicity,   genus
and reduction exponent of  I. }
\section{Introduction}
  \vspace{-0.3cm} Let ($R, \m, k$)  be a one-dimensional  local Noetherian
 domain   which is analytically irreducible and
residually rational.  In this paper we deal with the   {\it blowing-up  } $  \Lambda:= \LL(I)=
{\bigcup }_{n\geq 0}\ I^n:I^n  $ \   along a not principal \m-primary ideal  $ I$ of $R$.\\ The
problem of finding relations involving the multiplicity $e:=e(I)$,  the genus $\rho:=\rho(I)=l_R( \LL /R)$   and the
reduction exponent $\nu:=\nu(I)$,  was first studied    for $I=\m$  by Northcott  in the 1950s   and
later   by Matlis (see \cite{ma}), 
Kirby (see \cite{k}),  Lipman (see \cite{l}) and many others.      \par In this note we show that it is possible to describe
the   difference \ $2 \rho - e \nu$ in
terms of the {\it type sequence} $[r_1,....,r_n]$ of $R$  ($r_1$ is the {\it Cohen-Macaulay type}). Our main result
is the formula of Theorem \ref{primoint} in Section 4:
\vspace{-0.2cm}$$ \vspace{-0.3cm} 2 \rho = e \nu+\sum_{i\notin \Gamma} \ (r_i-1)-
d(R:\Lambda)-l_R(\LL^{**}/\LL)-l_R(R:\LL /I ^\nu). $$
Afterwards we use this statement   to improve classical results concerning  the equality
\vspace{-0.2cm}$$ \vspace{-0.2cm}
   R:\LL=I^{\nu} $$
which has been studied by several authors under the hypothesis that $R$ is Gorenstein.   Starting from a theorem of 
Matlis valid for
$\LL(\m)$ (\cite{m}, Theorem 13.4), Orecchia and Ramella  (\cite{or}, Theorem 2.6) proved that if the associated
graded ring
$G(\m)=\bigoplus_{n\geq 0}\m^n/\m^{n+1}$ is Gorenstein, then  $R:\LL=\m^{\nu}$.  Successively  Ooishi, in the case of 
the blowing-up along an ideal $I$,  
 proved that
$2\rho\leq
e\nu$ and that equality holds if and only if  $R:\LL=I^{\nu} $  (\cite{oo1}, Theorem 3). 
\par  In Section 5,  we consider the rings having  {\it type sequence} $[r_1,1,\dots,1]$ which are called almost
Gorenstein. For these rings we prove      that the   Ooishi's inequality $2\rho\leq e\nu\ \!$ 
becomes
$\ 2\rho\leq\nu e +r_1-1$ \  and that equality holds    if and only if   \   
$R:\LL=I^{\nu} $    (Theorem \ref{primointbis}).\\ \indent  In Section 6 we consider the case 
of the blowing-up along $\m$.   The study of the conductor
   $ \!  R:\LL $      provide some useful remarks when
\ $e =
\mu+1$  ($\mu$ is the embedding dimension of $R$)\ and when the  reduction exponent is 2 or 3. 

\section{Notations and Preliminaries.} \vspace{-0.3cm}
Throughout this paper    ($R, \m$)   denotes a one-dimensional local Noetherian
 domain with
residue field $k$. For simplicity, we assume that $k$ is an infinite field.
Let $\ \R\ $ be the
integral closure of
$R$ in its quotient field
$K $; we suppose that $\R$ is  a finite $R$-module and a DVR with a uniformizing parameter $t$,
which means that $R$ is analytically irreducible.  We also
suppose $R$ to be residually rational, i.e., \ $k
\simeq \R/t\R.$ \ We denote the usual    valuation
associated to $\R $
by\\
\centerline{$ v : K\longrightarrow \integ \cup
\infty,\qquad v(t)=1 .$}

\begin{(*)} \label{jager}  \ \\
 {\rm Under our hypotheses, for any fractional ideals $I
\supseteq J$ the
length of the $R$-module $I/J$ can be computed by means of valuations  (see
\cite{ma}, Proposition 1):}

 \centerline{ $l_R(I/J)=\# (v(I) \setminus v(J)).$ } \end{(*)}
Given two fractional ideals \ $I, J$ \ we define\   $I:J=\{x\in K \ | \ xJ \subseteq I\}$.

\begin{(*)} \label{defLL}   \ \\ {\rm
In the sequel we shall   consider  an \m-primary ideal  $I$ of $R$ which is not principal.   \vspace{-2mm}\\  The
{\it   Hilbert function   } and the  {\it Hilbert-Poincar\'e   series} of
$I$ are   respectively    \vspace{-0.2cm}
$$\vspace{-0.3cm} H_I(n)=  l_R (I^n/I^{n+1}), \ \  n\geq 0, \qquad P_I (z)=
\sum_{n\geq 0} H_I (n) z^n . $$ 
It is well-known that the power series \ $P_I(z)$
\ is rational: \vspace{-0.3cm}
$$\vspace{-0.3cm} P_I(z)=\displaystyle \frac {h_I(z)}{1-z},\quad \ {\rm
where}\quad \ h_I(z)= h_0+h_1z+h_2 z^2+....+h_{\nu}
z^{\nu}\in\integ[z], $$
  \ \  $h_0=l_R(R/I), \quad h_i= l_R(I^i/I^{i+1}) -
l_R(I^{i-1}/I^{i}),$  for all \ $i,\quad \ 1 \leq i \leq \nu$.
\vspace{0.1cm}\\
   The polynomial $h_I(z)$ is called the {\it h-polynomial} of $I$;\
moreover\par
  $ e(I):= h_I(1)$ is   the {\it multiplicity} of $I$, \par $ \rho(I):=
h_I'(1)$ is   called {\it genus} of
$I$,\ or   {\it reduction
number} of $R$\ \ if $I=\m$.  \vspace{0.1cm}\\  We shall say
that \
$h_I(z)$ is  {\it symmetric} \
  if  \   $h_i = h_{\nu - i}$ \
  for all \ $i, \ \ 0 \leq i \leq \nu$.  \vspace{2.5mm}\\
\noindent  The  {\it blowing-up  of \ $R$   along} \ $ I$ \      is
defined by  \vspace{-0.2cm}
$$\vspace{-0.3cm}\ \Lambda:= \LL(I)=
{\bigcup }_{n\geq 0}\ I^n:I^n    \qquad (cf.\cite{l}). $$
    Let $x\in I$ denote an element (called {\em a minimal
reduction of}   $I$) such that \\  $  I^{n+1}=xI^{n}$
\ for \ $n \gg 0$. \  Then
  (see \cite{l},  1):  \begin{enumerate}\item[(1)]   $x^n \LL =I^n \LL , \quad
\forall\ n
\geq 0$.

\item[(2)]
$e(I)=l_R(R/xR)=v(x)\geq H_I(n),$ for every $n\geq 0.$

\item[(3)] The   least   integer $\nu:= \nu(I)$\ such that $\ I^{n+1}=xI^{n}
\quad
\forall \ n\geq\nu, $ is called the {\it reduction exponent} of $I$.  It is
known that $ \ \nu(I)\leq e(I)-1$ and that the
following equalities hold: \\ $  \nu(I)\ =\  deg \ h_I(z) \ =\
min\big\{n\ |
\  l_R(I^n/I^{n+1})=e(I)  \big\}$ \par $  \ \ \ \ \ \ \  = min\{n \ | \  \LL =I^n:I^n\}\ =\
min\{n\ | \
I^n\LL=I^n\}$. 
 \item[(4)] $\rho(I) =l_R(\LL/R) $.\ \  Hence \vspace{-0.2cm}
$$\vspace{-0.3cm}l_R(R/I^n )=e(I) n-\rho(I),\qquad \forall\
n\geq\nu.$$
 
\item[(5)] If $h_I(z)$ is symmetric, then \    $l_R(R/I^{\nu})=
\displaystyle \frac{e(I) \ \nu}{2}$.  \\
This follows immediately from the fact that, if $h_I(z)$ is symmetric, then
\   $  2 \rho(I) = e(I) \ \! \nu $  \ (see the
proof of Lemma 3.3, \cite{oo3}).

\item[(6)] The inclusion \   $ R:\LL\supseteq I^{\nu}$ \ always holds and
the equality \ $R:\LL=I^n$
\ implies that \ $ n=\nu$ \ \ (\cite{oo1}, Proposition 1,
\cite{or}, Lemma 1.5).
\end{enumerate}}
\end{(*)}

\begin{(*)} \ \\ {\rm We shall consider also:}  \par 
\vspace{ 0.2cm}     
\  \ $\begin{array}{rcl} v(R) &:=&  \{v(x),\;
x\in R,\; x\neq 0\} \subseteq \nat ,  {\rm    \ the} \  numerical\
semigroup\ {\rm of} \   R. \\
 \gamma_R&:=&R:
\R,  \ {\rm   \ the {\it \ conductor\ ideal }\ of }\ R.\\
 c&:=& l_R(\R/\gamma_R), \ {\rm  \  the} \ conductor \ {\rm of} \
 v(R) , {\rm \ \ such\ that\ }
 \ \gamma_R  =t^c\overline R .\\
 \delta&:=& l_R(\R/R),   {\rm  \ the\ }\ singularity\
degree\ of\ R.\\
n&:=&c-\delta=l_R(R/\gamma_R). \end{array}\vspace{- 0.2cm}    $ 
\end{(*)}

\begin{(*)}\label{modcan}   \ \\{\rm  In our hypotheses $R$ has
a  {\em canonical module}  \  $\w$,    unique up to isomorphism.\\
We list below some well-known properties of  $\w$,   useful in the sequel
(see \cite{hk}).     We always assume   that \
$R\subseteq
\omega \subset \R.$ \begin{enumerate}

\item[(1)] $\w:\w=R$ \ and \  $\omega:(\omega:I)=I$ \ for every fractional
ideal  $I $.

\item[(2)] If \  $I \supseteq  J$, \ then \ $l_R(I/J) = l_R(\w:J/\w:I).$

\item[(3)] $ v(\w) = \{ j\in \integ \enskip | \enskip c -1
 -j \notin v(R) \}  $, \ hence \ $c-1 \notin v(\w)$ \ and \ $c+\nat \subseteq v(\w)$.

\item[(4)] $R$ is Gorenstein \  {\it if and only if} \ $    \w  =  R$ \
{\it  if and only if} \  $ R: \w = R.$   \\
Otherwise\qquad $ \gamma_R \subseteq R: \w \subseteq \m.$
\item[(5)]   (see \cite{ooz}, Lemma 2.3). For every  fractional ideal $I$, \\
\vspace{-0.3cm}
$$\vspace{-0.2cm} s\in v(I\w)  \quad  if \ and\ only\ if \quad   c-1-s\notin v(R:I).$$
\end{enumerate}}
\end{(*)}

\begin{(*)} \label{ts}  \ \\
{\rm We   recall the notion of {\em type sequence} given for
rings by Matsuoka in 1971, recently revisited in \cite{bdf} and extended
to modules in \cite{oz1}.\\
Let $n:=c-\delta ,$ and let $s_0=0<s_1<\ldots <s_n=c$ \ be the first
$n+1$ elements of $v(R).$ \ For each $i=1,\dots,n$, \ define the ideal
$R_i:=\{ x\in R\; :\;
v(x)\geq s_i\} $ \ and consider the chains:
\vspace{-0.2cm}  $$R=R_0\supset R_1=\m\supset R_2\supset
\ldots \supset R_n =\gamma_R \vspace{-0.2cm}$$
$$R\; =\; R\; :\; R_0\; \subset R\; :\; \m\subset
\; R\; :R_2\; \subset \; \ldots \subset R\; :R_n=\; \R$$
 For every $i=1,\ldots, n,$ \ put \
 $r_i:=l_R(R:R_{i}/R:R_{i-1})=l_R(\w R_{i-1}/\w R_{i}).$\par

\noindent The {\em type sequence} of $R,\; $ denoted by $\; t.s.(R),$ is
the sequence
$\; [r_1,\dots,r_n].$ \par

\noindent We list some properties of type sequences  useful in the sequel
(see
\cite{bdf}): \begin{enumerate}

\vspace{-0.2cm} \item[(1)] $r:=r_1$ is the {\em Cohen-Macaulay type} of $R$.

\vspace{-0.2cm} \item[(2)] For every $\ i=1,\dots,n,$ \ we have \
$ 1\leq r_i \leq  r_1 $.

\vspace{-0.2cm} \item[(3)] $  \delta   = \sum_1^n \  r_i $,  \ \ and \ \ \ $2
\delta -c = l_R(\w/R)=
\sum_1^n \
(r_i-1)$.

\vspace{-0.2cm} \item[(4)] If \ $s_i \in v(R: \w),$ then
the correspondent $\ r_{i+1}$ is $1 $ (see \cite{ooz}, Prop.3.4).
\end{enumerate}
}\end{(*)}

\begin{(*)} \ \\
{\rm  We recall that  ring $R$ is called  {\it almost Gorenstein} if it
satisfies the equivalent conditions\par
(1) $ \ \m =   \m\ \!\w.$ \par
(2) $\ r_1-1= 2 \delta - c. $ \par
(3) \ $ R: \w \supseteq \m.$\\
By the above property \ref{ts},(3), it is clear that $R$ is   {\it almost
Gorenstein} if and only if  $ t.s.(R)=
[r_1,1,\dots,1] $  and that {\em Gorenstein}  means  {\em almost
Gorenstein} with $r_1=1 $. }
\end{(*)}

\begin{(*)}\label{riflessivi}   \ \\{
\rm
For any fractional ideal $I$ of $R$ we set \ $I^*:=R:I $. \
 Notice that:\vspace{-0.2cm}
$$  \vspace{-0.2cm} I \subseteq I^{**}  \subseteq   I\w.$$   In fact, \ $I^{**} =
R:(R:I) \subseteq \w :(R:I) = I\w.$ 
     }\end{(*)}

\begin{(*)}\label{chiusure}   \ \\{\rm
We recall that  the  {\it integral closure} of   an ideal \ $I$ of $R$   is   \
$ \Ib := I \R \cap R$  and that  
 $I$ is said to be
 {\it integrally closed } if   \ $I=\Ib$.\\ In \cite{oo2} Ooishi characterizes curve singularities which can be
normalized by the first blowing-up  along the ideal $I$ in terms of integral
closures:\vspace{-0.2cm}
$$ (*)\qquad\LL   = \R \quad if \ and \ only \ if \quad I^n = \overline{I^n} \quad
{\rm for\ all }\ n \geq \nu.\vspace{-0.2cm}$$
We introduce a weaker notion of closure, namely the  {\it canonical closure} of   $I$  as \ $ \It := I\w \cap R$. \ We'll see that this   notion is
particularly meaningful for almost Gorenstein rings.
  \ Recalling \ref{riflessivi}, we can easily see that  \ $I \subseteq \it {I}^{**}\subseteq\It
\subseteq \Ib$, \ so \vspace{-0.2cm}$$\vspace{-0.2cm}I=\Ib \quad  implies\  that  \quad  I=\it {I}^{**}=\It.$$ \  For the
canonical closure the analogue of statement   $ (*)$  is: 
\vspace{-0.2cm}
 $$ \quad \LL   = \w \LL \quad if \ and \ only \ if \quad I^n= \widetilde{I^n}
\quad  {\rm for\ all }\ n \geq \nu.$$
  This fact is shown in the next proposition.}
\end{(*)}    
\begin{prop} \label{scopp.riflessivo}\  \\
Let \
$\LL:=\LL(I)$ be as above.
  We have the following   groups of equivalent conditions:
\begin{enumerate}
\item[ {\rm (A)}]    \quad ${\rm(A_1) }\quad   \w \subseteq \LL;  \
\ \quad{\rm(A_2) }\quad \  \w
\Lambda =
\Lambda; \ \ \quad{\rm(A_3) }\quad \  \w: \Lambda =
R: \Lambda$;
\item[] \quad ${\rm(A_4) }\quad \ I^n=\widetilde{I^n}\quad \forall \ n \geq \nu;  \ \hspace{1.5cm}
  {\rm(A_5) }\quad \ \w I^n =  I^n \quad \forall \ n \geq
\nu$;
  \item[] \quad ${\rm(A_6) }\quad there \ exists \  n > 0 \   such \ that \   \w  I^n=  I^n$.
\item[{\rm (B)}]   \quad $ {\rm(B_1) }\quad   \LL = \LL^{**};  \   \qquad
{\rm(B_2) } \quad  \  \w:\LL = \w(R:\LL).$
\end{enumerate}
Moreover the following facts are equivalent\begin{enumerate} \item[(1)] \  Conditions {\rm (A)} hold.  \item[(2)] \
Conditions {\rm (B)} hold and \
$ R:\LL\subseteq R:\w$.
\end{enumerate}
\end{prop}
{\it Proof} \\ Let's begin to prove that the equalities
 $  \ I^n=  \widetilde{I^n } \quad \forall \ n \geq \nu$ \ imply that \ $
\w \subseteq \LL$. \\
 Let $k_1$ be the minimal exponent such that \ $ I^{k_1}
\subseteq R: \w$ \ ($k_1$ exists since $R: \w \supseteq \gamma_R$).
\  If \ $k_1
\geq \nu$, \ then   \  $I^{k_1}=  \widetilde{I^{k_1} }= \w
I^{k_1}  \cap R = \w  I^{k_1}$ \ and this yields  \ $ \w \subseteq
I^{k_1}:I^{k_1}= I^{\nu}:I^{\nu} = \LL$. \ If \ $k_1 < \nu$, \ then \ $
I^{\nu} \subseteq I^{k_!} \subseteq R: \w$, \ hence \ $  \w
I^{\nu}
\subseteq R$. \ Thus, \ $ I^{\nu}= \widetilde{I^{\nu} }=  I^{\nu}\w  \cap R
= I^{\nu}\w$, \ which means \ $ \w \subseteq  \LL$. \\
All the other implications in group (A) and also that ones  in group (B) hold
by the properties    of the canonical module.\\
To prove {\it (A) implies (B)}, note that \ by \ref{riflessivi} \ $\LL^{**} \subseteq  
  \w\LL = \LL$.\\
Moreover, if (A) holds, then $(R: \LL) \w \subseteq (R: \LL) \LL \subseteq
R$, \ hence \ $R: \LL
\subseteq R: \w$.\\
     Under the further assumption   \ $R: \LL
\subseteq
R: \w$,\ we can prove {\it (B) implies (A)} because the fact
\
$\w:\LL =\w(R:\LL)\subseteq
\w(R: \w)\subseteq R$ \ leads to \ $\LL\supseteq\w$.  

\begin{rem} \label{scop} \end{rem} \begin{enumerate}\item[(1)] If $R$ is almost Gorenstein, then \ $R:\LL\subseteq R:\w$,
hence conditions (A) and (B) above are equivalent.
\item[(2)]If $I$ is a canonical ideal, i.e., $I \simeq \w$, then
conditions {\rm (A)} and
{\rm (B)} hold, because  $\LL$ is
reflexive and  \ $R: \LL\subseteq R: \w$ \ (see  \cite{ooz}, Remark 2.5).\end{enumerate}

\section{The first Formula.}\vspace{-0.3cm}
 In the following we use  the   notation  introduced in Section 2.\\
 $\Lambda:= \LL(I)= \bigcup _{n>0}  I^n:I^n$ is  the {\it blowing-up of
$R$} in an
$\m$-primary ideal   $I$ which is not principal and  \ $e:=e(I), \  \nu:= \nu(I), \ \rho:= \rho(I)$ \
are respectively the
{\it multiplicity}, the {\it reduction exponent} and the {\it genus}   of
$I$. \\
   Moreover we consider \ \  $\gamma_R:=R   : \R, \qquad  \delta: = l_R(\R/R), 
\qquad \     c:=l_R(\R/\gamma_R),$
\par \hspace{3cm}
$ \gamma_{\LL}:= \LL:\R,\qquad\delta_{\LL}:=l_R(\R/\LL),\qquad  c_{\LL}:=  l_R(\R/\gamma_{\LL}) 
$.\\
Finally,   $x\in I$  denotes a {\it minimal reduction} of $I$.

\begin{(*)} \label{notations} {\rm We begin with a few remarks involving  the conductor
ideals respect to the canonical inclusions \ $R \subseteq \LL \subseteq\R$. \ We have the following diagram:}
$$\begin{array} {ccc} \gamma_{\Lambda} & \  & \ \\
\cup\hspace{.3mm}\rule{.1mm}{2.1mm}   &\ & \\
\gamma_R&\quad\subseteq\quad & R:\Lambda \ \  \\
\cup\hspace{.3mm}\rule{.1mm}{2.1mm}&&\cup\hspace{.3mm}\rule{.1mm}{2.1mm}\\
(R:\Lambda) \ \gamma_{\Lambda}&& I^{\nu}\\
\cup\hspace{.3mm}\rule{.1mm}{2.1mm}&&\cup\hspace{.3mm}\rule{.1mm}{2.1mm}\\
x^{\nu} \ \gamma_{\Lambda}&=&I^{\nu} :\R
\end{array}$$
\end{(*)}
 
\begin{prop} \label{conduttori}  \
\begin{enumerate}
\item[(1)]  \ $ c - c_{\LL}\leq e  \nu$.
\item[(2)]  \ $l_R(R/\gamma_{R})-l_R(\LL/\gamma_{\LL}) =  c - c_{\LL}-
\rho=e\nu-\rho-l_R(\gamma_R/x^{\nu}\gamma_{\Lambda}) \leq l_R
(R/I^{\nu})$.
\item[(3)]   \ The following facts are equivalent:
\begin{enumerate}
\item $c -  c_{\Lambda}= e  \nu. $
\item $\gamma_R = x^{\nu}\gamma_{\Lambda}.$
\item $\gamma_R \subseteq I^{\nu}.$
\item $l_R(R/\gamma_R)-l_R(\LL/\gamma_{\Lambda})=l_R(R/I^{\nu}).$
\end{enumerate}
\end{enumerate}
\end{prop}
{\it Proof} \begin{enumerate}
\item[(1)]  Considering the diagram in \ref{notations} we see that:\\
 $ c-c_{\Lambda}=l_R(
\gamma_{\Lambda}/\gamma_R)=l_R(\gamma_{\Lambda}/x^{\nu}\gamma_{\Lambda})-l_R(\gamma_R/x^{\nu}\gamma_{\Lambda})= e\nu
-l_R(\gamma_R/x^{\nu}\gamma_{\Lambda}).$
 \item[ (2)] Since \ $\rho= \delta - \delta_{\Lambda}, $ \
using part (1) of the proof       we obtain: \par
 \ \ \ \ \ \ $l_R(R/\gamma_R)-l_R(\Lambda/\gamma_{\Lambda}) = (c-  \delta)
-(c_{\Lambda}-\delta_{\Lambda}) = c -
c_{\Lambda}-
\rho $ \par \ \ \ \ \ \ \ \ \ \ \ \ \ \ \ \ \ \ \ \ \ \ \ \ \ \ \ \ \ \ \ \
\ \ \ $  =   e\nu -\rho
-l_R(\gamma_R/x^{\nu}\gamma_{\Lambda})\leq l_R(R/I^{\nu}).$
 \item[ (3)] Equivalences \
{\it (a) if and only if  (b)  }\ and \ {\it (b) if and only if  (d)} \ are
immediate by item (2).\\
  To prove {\it  (b) implies (c)}, we note that  \  $\gamma_R  =
x^{\nu}\gamma_{\Lambda}=I^{\nu}:\R\subseteq   I^{\nu}$.
  Conversely, assumption (c) implies that \  $\gamma_R  =\gamma_R:\R
\subseteq I^{\nu}:\R \subseteq \gamma_R$, \   hence
$\ \gamma_R \ =\ I^{\nu}:\R \ =\ x^{\nu}\gamma_{\Lambda}.$ \end{enumerate}

\begin{rem}  \label{nonaltro} \
\begin{enumerate}  {\rm
\item[(1)] In view of item (2) of the above proposition   we have the
inequality   \vspace{-0.2cm}
$$l_R(R/\gamma_{R})-l_R(\LL/\gamma_{\LL}) \geq -\rho\vspace{-0.2cm}$$
    and, in the case  $I=\m$, \ $
l_R(R/\gamma_{R})-l_R(\LL/\gamma_{\LL})\geq e-\rho $.\\ In    Example
\ref{negativo} we show
that both   these minimal values can be reached.
\item[(2)] Conditions (3) of \ref{conduttori}  imply that
$R:I^{\nu}\subseteq \R$,
but if this inclusion holds we need not
have the above  equivalent conditions (see   Example \ref{cinque}).
\item[(3)] Conditions (3) of \ref{conduttori} \ imply the Conductors
Transitivity
Formula:\vspace{-0.3cm}
$$\ \gamma_R =(R:\Lambda) \
\gamma_{\Lambda}.\vspace{-0.3cm} $$  Example \ \ref{ninon2}   shows that the
converse does
not hold.
\item[(4)]  Conditions (3) of \ref{conduttori} do
not imply that
$ R:\Lambda=I^{\nu}$. \ This can be seen in Example \ref{anegativo};
however   next
lemma   shows that  the
converse is true. } \end{enumerate}
\end{rem}

 \begin{lemma}  \label{conduttore=potenza} \ \\ If \
$R:\LL=I^{\nu}$, \ then we have
\begin{enumerate}

\item[(1)] \ The  equivalent conditions   of Proposition
\ref{conduttori},(3)  hold.
\item[(2)] $ \ \LL^{**}=\LL $.
\end{enumerate}
\end{lemma}
{\it Proof} \\
(1)  It is  clear considering the diagram in \ref{notations}. To prove part 
(2), observe that  condition  \ \ $ R:\LL=I^\nu =  \LL I^\nu$ \ implies   $ \
\LL^{**}= R:\LL I^{\nu}=I^\nu:I^\nu=\LL$.    \vspace{0.3cm}  \\ 
 From the above considerations we obtain   a first formula connecting the
invariants $\rho,e, \nu$ associated to the
ideal $I$ with the invariants $c,  \delta  $ of  $R$  by means of  the length of the quotient $R:\LL/I^{\nu}$. 
This formula will
be successively improved in Theorem
\ref{primoint} by using type sequences.

\begin{prop} \label{modcanLL} \
\begin{enumerate}
\item [(1)] \
 $  2 \rho= e \nu + (2 \delta - c)-l_R(R:\LL/I^{\nu}) - l_R(\w \LL/\LL)$.
\item [(2)] \    The following facts are equivalent:\vspace{-0.1cm}
\begin{enumerate}
\item \ $2 \rho = e  \nu +(2 \delta - c)$. 
\item  \ $\LL$   is Gorenstein and \ $c - c_{\LL}= e  \nu$.
\item  \ $R:\LL=I^{\nu}$ \ and \ $\w \LL= \LL$.
\item  \ $R:\LL=I^{\nu}\subseteq  R:\w.$
\end{enumerate}
\end{enumerate}
\end{prop}
{\it Proof} \\
  From \  $2 \rho= 2 \delta - 2 \delta_{\LL}=  2 \delta - c -(2
\delta_{\LL}- c_{\LL})+c- c_{\LL}+e \nu -e \nu$,   we get  \vspace{-0.1cm}
$$(*) \qquad  2 \rho= e \nu + (2 \delta - c)-(2 \delta_{\LL}-
c_{\LL}) - (e \nu -c +c_{\LL})\vspace{-0.1cm}$$
Hence  the equivalence \ {\it  (a) if and only if (b)} \ of (2) is clear.\\
Since \  $I^{\nu} \subseteq R:  \LL  \subseteq R,$ \ we have  \vspace{-0.1cm}
$$ (**) \ \ \ \ \ \ \ \
l_R(R/R:\LL)=l_R(R/I^{\nu})-l_R(R:\LL/I^{\nu}) = e
\nu -
\rho -l_R(R:\LL/I^{\nu}).\vspace{-0.2cm}$$
From the inclusions  \ $R \subseteq   \LL  \subseteq \w \LL$ \ and  \
$R \subseteq \w \subseteq \w \LL$, \ we obtain that \vspace{-0.1cm}
 $$ l_R(R/R:\LL)=  l_R(\w \LL/\w) =l_R(\w \LL/\LL)+ \rho - (2 \delta -
c).\vspace{-0.2cm}$$
Substituting this in the first member of $(**)$ we get the first formula and
also the equivalence {\it (a) if and only if (c)}.\\
Finally, {\it (c) if and only if (d)} follows by  using Proposition \ref{scopp.riflessivo}.

\section{Formulas involving type sequences.}\vspace{-0.3cm}
We keep the   notation    of the above section.
We have   seen in \ref{modcanLL} that \vspace{-0.2cm} $$2 \rho  \leq
e   \nu  +(2 \delta - c).\vspace{-0.2cm} $$
Using the  notion of   {\it type sequence} we    insert a new term
in this inequality  (see Theorem   \ref{primoint}):
\vspace{-0.2cm} $$2 \rho
\leq  e
\nu +\sum_{i\notin
\Gamma } \ (r_i-1)  \leq
e   \nu  +(2 \delta - c).\vspace{-0.2cm} $$
 We study also conditions to have equalities. To do this
we  introduce   the positive invariant  $ d(R:\LL)$,
which plays a crucial role in this context.

\begin{defin} \ \\
Let, as above,  $s_0=0,s_1,\ldots ,s_n=c$ be the first
$n+1$ elements of $v(R)$,   $ n=c-\delta. $ \  Let \
$t.s.(R)=[r_1,...,r_n]$ \ be the type sequence
of $R$.   \
We call $\ d(R:\LL)\ $ the number
\vspace{-0.2cm}
$$ d(R:\LL):= l_R(\R/\LL^{**}) -\sum_{i\in \Gamma}r_i \vspace{-0.2cm}$$
where $\Gamma$ denotes the numerical set
\ \ $\Gamma := \{ i\in \{1,..,n\} \ | \ s_{i-1} \in v(R:\LL)\}$.\end{defin}

\noindent Note that \vspace{-0.2cm}
$$ \# \Gamma=l_R(R:\LL/\gamma) = l_R(\R/ \w \LL)   $$

The following proposition ensures  that \ $ d(R:\LL) \geq 0$.

\begin{prop} \label{diseg}   \  \\
 We have 
$$ l_R(\R/\w \LL)\leq \sum_{i\in \Gamma}\ r_i  \leq
l_R(\R/ \LL^{**}) .  $$
 \end{prop}
{\it Proof} \\   The first inequality is obvious  since $r_i \geq 1 \ \
\forall \ i$. For the second one  we shall use
  property (5)  of \ref{modcan}   with $I =R:\LL$: \\
 \centerline{$s\in v(I\w)$  \ {\it if and only if} \ $ c-1-s\notin
v(\LL^{**}) $.}
If $x_{i-1}\in I $ is such that $ v(x_{i-1})= s_{i-1},$ then by definition
\vspace{-0.2cm}
$$r_i= l_R (\w R_{i-1}/\w R_i) =
l_R(x_{i-1}\w+\w R_{i}/\w R_{i})=\#
\{v(x_{i-1}\w +\w R_{i})\setminus v(\w R_i)\}.\vspace{-0.2cm}$$   Since $
\ v(x_{i-1}\w)
\subseteq  v(I\w)$,  \  the assignment
\
$y
\to c-1- y$
\ defines an injective map
\vspace{-0.2cm}
$$  \vspace{-0.2cm}\bigcup_{i \in \Gamma }\{ v(x_{i-1}\w +\w R_{i})\setminus v(\w R_{i}) \}
\longrightarrow \nat \setminus v(\LL^{**}).$$
 From the fact that the numerical
sets  \vspace{-0.2cm}
$$  \vspace{-0.2cm}\{ v(x_{i-1}\w +\w R_{i})\setminus v(\w
R_{i})\}, \ \ i\in \{1,\dots,n\}, $$ are disjoint by construction   we deduce that
\vspace{-0.2cm} $$\sum_{i\in \Gamma }  r_i \leq l_R(\R/\LL^{**}).
\vspace{-0.2cm}$$
The next proposition collects some useful properties of the invariant
$d(R:\LL) $ and allows us to find  sufficient
conditions to have
$d(R:\LL)=0$.

\begin{prop} \label{dopodefc}   \ \\
Let  $i_o\in\nat$ be such
that \
$e(R:\LL)=s_{i_0}$. Then
\begin{enumerate}
\item[(1)] \    $ d(R:\LL)=l_R(\w \LL / \LL^{**} )-   \displaystyle{ \sum_{i\in
\Gamma }\ ( r_i-1) }$.
\item[(2)] \  If \ $\w \subseteq \LL^{**}$,  i.e.,  \ $R:\LL  \subseteq R:\w$, \ then \ $ d(R:\LL)=0$.
\item[(3)]  \ $  d(R:\LL)= \displaystyle{\sum_{i> i_0, \ i\notin \Gamma} }  r_i\
-\ l_R(\LL^{**}/R_{i_0}^*)$.
\item[(4)]  \ If \ $R:\LL$ is integrally closed, then \ $d(R:\LL)=0$.
 \end{enumerate}
 \end{prop}
{\it Proof} \begin{enumerate}
 \item[(1)]   $ d(R:\LL) = l_R(\R/\LL^{**}) -\displaystyle \sum_{i\in
\Gamma}r_i =l_R(\w \LL / \LL^{**}
)-(\sum_{i\in \Gamma}r_i -l_R(\R/\w\LL) ) $.
\item[(2)] The inclusion \ $\w \subseteq \LL^{**} $   implies that \ $\w \LL =\LL^{**}$,   hence  the thesis by (1), recalling that $d(R:\LL)\geq 0$.
\item[ (3)]  After writing \  $l_R(\R/\LL^{**})=
l_R(\R/R_{i_0}^*)-l_R(\LL^{**}/R_{i_0}^*),  $     the thesis is
clear  since \
\vspace{ -0.1cm} $$l_R(\R/R_{i_0}^*)=\displaystyle  \sum_{i> i_0}   r_i.
\vspace{ -0.3cm}$$
 \item[(4)]  This results from the above item, because the fact that $R:\LL$ is
integrally \par closed means that \ $R:\LL= 
R_{i_0}$.\end{enumerate}
The next theorem provides a link between the {\it type sequence} of $R$ and
the {\it genus}  $\rho $  of the
ideal
$I$.
\begin{thm} \label{bounds} \ \begin{enumerate}
\item[(1)] \  $\rho =   \displaystyle  \sum_{i\notin
\Gamma}r_i-l_R(\LL^{**}/\LL)-d(R:\Lambda) \leq r \ l_R(R/R:\LL)$.
\item[(2)] \   Let  $i_o\in\nat$ be such
that \
$e(R:\LL)=s_{i_0}$. Then 
\item[] $\ \rho=\displaystyle \sum_{i\leq i_0} r_i  -l_R(\LL^{**}/\LL)  +
l_R(\LL^{**}/R_{i_0}^*) $. \end{enumerate}
 \end{thm}
{\it Proof} \begin{enumerate} \item[(1)]    From the inclusions \ $R
\subseteq \LL
\subseteq \LL^{**} \subseteq \R$ \ we obtain \par \ \ \ \
$\rho=l_R(\LL/R)=\delta -
l_R(\LL^{**}/\LL)- l_R(\R/\LL^{**})$ \par $  \ \ \ \ \ \ \ = \delta -
l_R(\LL^{**}/\LL)- d(R:\Lambda)-  \displaystyle  \sum_{i\in \Gamma}r_i
$. \\ Thus the first equality is clear since  \ $\delta -
  \displaystyle  \sum_{i\in \Gamma}r_i= \displaystyle  \sum_{i\notin
\Gamma}r_i   $. \par  The   inequality follows  immediately, recalling that
\ $r_i\leq r \ \ \forall \ i $ \ and that
\vspace{-0.2cm}$$\vspace{-0.2cm} l_R(R/R:\LL)=\#\big(\{1,\dots ,n\}\setminus\Gamma \big).$$  
\item[(2)]   By substituting  formula (3) of   \ref{dopodefc} in
   formula (1) above, we obtain \par
 \ \ \ \ $\rho =  \displaystyle \sum_{i\notin \Gamma}  r_i
-l_R(\LL^{**}/\LL)-\sum_{i > i_0,  i\notin \Gamma }   r_i  +
l_R(\LL^{**}/R_{i_0}^*)=$ \par $ \ \ \ \ \ \ =
 \displaystyle \sum_{i\leq i_0} r_i  -l_R(\LL^{**}/\LL)  +
l_R(\LL^{**}/R_{i_0}^*) $.  \end{enumerate}

\begin{rem} \ \\
{\rm   In the case   $\LL=\R$ \   the inequality  \  $\rho   \leq r \
l_R(R/R:\LL)$  of   Theorem \ref{bounds}
   \ gives
\ the well-known relation \ $\delta \leq r \ (c- \delta)$ \ (\cite{ma},
Theorem 2).\\
The maximal value   $\rho  = r \ l_R(R/R:\LL)$   is achieved if and only if \
$  r_i=r$ \
for all $
\ i \notin \Gamma, \  \LL= \LL^{**} \ $ and $\ d(R:\LL)=0$; \  this
    happens for instance
if \  $I=\m$ \ and \ $e=\mu$ (see \ref{e=mu}), \ or if \
   $R$ is   {\it Gorenstein}. }
\end{rem}

\begin{coro} \label{minor} \
\begin{enumerate}
\item[(1)]  \ $ e \nu  +rl_R(R:\LL/ I^{\nu}) \leq (r+1)\ l_R(R/I^{\nu})$.
\item[(2)] \ If the  h-polynomial is symmetric, then \ $   l_R(R:\LL/
I^{\nu})\leq \displaystyle
\frac{ r-1 }{r}
\cdot \displaystyle \frac{e \nu }{2}$
\end{enumerate}
\end{coro}
{\it Proof} \begin{enumerate}
\item[(1)]  From the first item   of the theorem   we have: \

$   \rho = e\nu - l_R(R/I^{\nu})   \leq r l_R(R/R:\LL)=rl_R(R /I^{\nu})- r
l_R(R:\LL/I^{\nu}).$

  The thesis follows.
\item[(2)] By   property (5) of \ref{defLL} it suffices to substitute \
$l_R(R/I^{\nu}) =
\displaystyle \frac{e \nu }{2}$ \ in (1).\end{enumerate}

\begin{thm} \label{primoint} \ \begin{enumerate} \item[ (1)] \  $2 \rho = e
\nu+\sum_{i\notin \Gamma} \ (r_i-1)-
d(R:\Lambda)-l_R(\LL^{**}/\LL)-l_R(R:\LL /I ^\nu).$
\item[(2)]  \    The following facts are equivalent:\vspace{-0.1cm}
\begin{enumerate}
\item \ $  2 \rho = e
\nu+\sum_{i\notin \Gamma} \ (r_i-1).$ \item \ $R:\LL =  I ^\nu
$    and   \ $d(R:\Lambda)=0$.
\end{enumerate}\end{enumerate}
\end{thm}
{\it Proof} \\  (1). \  We can rewrite   formula (1)  of Proposition \ref{modcanLL}
as:\vspace{-0.2cm}$$\vspace{-0.2cm}    2 \rho =  e \nu+ \displaystyle \sum_{i\notin \Gamma} \ (r_i-1)+\sum_{i\in\Gamma} \
(r_i-1) - l_R(R:\LL /I ^\nu)-l_R(\LL^{**}/\LL)- l_R(\w\LL/\LL^{**}).$$
 So using    item (1) of Proposition \ref{dopodefc}, we
obtain part (1).  \\  (2) follows
from part (1) by virtue of  Lemma
\ref{conduttore=potenza} recalling that $d(R:\LL)\geq 0 $.\vspace{0.3cm}  \par
 We remark that the equality \ $R:\LL =  I ^\nu$ \
 does not ensure that \ $d(R:\Lambda)=0$ \ (see Example \ref{cnonnullo}).\par

\section{Almost Gorenstein Rings.}\vspace{-0.3cm}

In this section we deal with {\it almost Gorenstein} rings.  The notations
will be the same as in  the preceding sections.\par
Under the hypothesis {\it R almost Gorenstein},
  the  formulas  in   
\ref{modcanLL}, \ref{bounds} and \ref{primoint} involving the genus $\rho(I)$    
are considerably simplified and allow us to extend some well-known results
concerning the equality   $     R:\LL=I^{\nu}.$ 
  Recently Barucci and Fr\"oberg stated the
equivalence    \ {\it (a) if and only if (c)} \
of    next Theorem\
\ref{primointbis} in the
 case {\it R almost Gorenstein} and \ $\LL=\LL(\m)$ \  (see \cite{bf},
Proposition 26). \vspace{0.3cm}\\
First, inspired by the famous   result of Bass
"{\it A one-dimensional Noetherian local domain $R$ is Gorenstein
if and only if each nonzero fractional ideal of $R$ is reflexive}" (see
\cite{ba}, Theorem 6.3),  \  we
notice that:

\begin{prop} \label{Matlis2}    \ \\
 $R$ is  almost  Gorenstein \  if and only if  \ $
  \w J =J^{**}$  for  every not  principal fractional ideal  $J$. \end{prop}
{\it Proof} \\
Suppose  $R$   almost  Gorenstein. By \ref{riflessivi} it suffices to prove
that \ $
\w J \subseteq J^{**}$. \ Since \ $R:J=\m:J$, \ we
have
\ $ (R:J)J
\w
\subseteq \m \w =\m$, \ hence \ $  \w J \subseteq J^{**}$. \\
 The opposite implication follows immediately by taking  $J= \m$.

\begin{coro} \label{overring}   \ \\If  $R$ is an almost Gorenstein ring, \ 
then \begin{enumerate}
\item[(1)]
  $ \LL^{**}= \w \LL$ \  and \  $ d(R:\LL)=0$.
\item[(2)] $\rho =  r-1+l_R(R/R:\LL) - l_R(\LL^{**}/\LL)$.
\end{enumerate}
\end{coro}
{\it Proof} \begin{enumerate}
\item[(1)]    The second equality
follows from  Proposition        \ref{dopodefc}, (1).
\item[(2)] Apply Formula (1) of Theorem \ref{bounds}, observing that
in the almost   Gorenstein case \vspace{-0.2cm}$$ \vspace{-0.2cm} 
\displaystyle
\sum_{i\notin
\Gamma}r_i=  r-1+l_R(R/R:\LL) .$$ \end{enumerate}

\noindent Under the  assumption $R$ {\it almost Gorenstein},  since \vspace{-0.2cm}$$
\vspace{-0.3cm}\w\LL=\LL^{**},\quad
\sum_{i\notin\Gamma} (r_i-1)=r-1=2\delta-c  \quad{\rm and}\quad 
   d(R:\LL)=0 $$  both  Proposition \ref{modcanLL} and Theorem
\ref{primoint} give the next theorem.
\begin{thm} \label{primointbis}   \ \\
Assume that
$R$ is an almost Gorenstein ring and let \  $\LL =\LL(I)$. 
Then:
\begin{enumerate}
\item[(1)]  $2 \rho = e  \nu+r-1-l_R(R:\LL / I^{\nu})-l_R(\LL^{**}/\LL). $
\item[(2)]  The following conditions are equivalent:
\begin{enumerate}
\item
$2 \rho = e  \nu+r-1 $.
\item  $\LL$ is Gorenstein and \ $c -  c_{\LL}= e  \nu$.
\item $ R:\LL = I^{\nu}$.\item  $ \w:\LL = I^{\nu}$.

 \end{enumerate}
 In this case   the equivalent conditions (A)  of   Proposition \ref{scopp.riflessivo} hold. 
\end{enumerate}
\end{thm}
{\it Proof} \ \\
We have only  to prove   {\it (c) if and only if (d)}.\\      {\it (c) implies (d)}.
 \ By Lemma  \ref{conduttore=potenza}  we have \ $\LL=\LL^{**}=\w\LL $. \   Hence 
  \ $\w:\LL=R:\LL$ by duality.   
 \\ To prove {\it  (d) implies (c)}, we notice that   $I^{\nu} \subseteq R: \LL
\subseteq\w:\LL $. 
 \begin{coro} \label{minorbis}
 \ \\
If
$R$ is an almost Gorenstein ring   and the h-polynomial is symmetric, then \par
 $  l_R(R:\LL/I^{\nu}) \leq r-1$ \ \  and \ the equality holds \ 
  if and only if   \ \  $\LL=\LL^{**}$. 
\end{coro}
{\it Proof} \\
   The symmetry of the {\it h-}polynomial gives  \ $2\rho = e \nu  $   (see \ref{defLL}(5)), hence it suffices to substitute
this in formula (1) of the theorem.\\ \vspace{-1mm} 

 We note that   the condition \ $
l_R(R:\LL/I^{\nu})= r-1$ \  does not imply that the {\it h-}polynomial is symmetric: see for instance Example
\ref{comeGnonGor},  where $R$ is {\it almost Gorenstein} with
$r(R)>1$ and   Example
\ref{GnonGor}, where
$R$ is Gorenstein.
  Example \ref{GnonGor}   shows also that the hypotheses \ $R$ Gorenstein and
\ $  2\rho=e\nu$ \    do  not
give the symmetry of the $h$-polynomial.\par \vspace{2.5mm}
 The following statement of Ooishi (see \cite{oo1}, Corollary 6) can be
obtained as a direct consequence of our preceding
results.

\begin{coro}\label{Ggor} \ \\
 If $R$ is Gorenstein and the h-polynomial is symmetric, then the equivalent
conditions  (2) of Theorem \ref{primointbis} hold.
\end{coro}\par
   Another immediate consequence of Theorem \ref{primointbis} is the natural
generalization of Theorem 10 of \cite{oo1} to the
almost Gorenstein case.

\begin{coro} \label{sottocaso}   \ \\ Suppose $R$ almost Gorenstein.\\
The equality \
 $\gamma= I^{\nu}$ \ holds \ \ if and only if \ \ $\LL=\R$ \   and  \  $2 \delta= e
\nu +r -1$.
\end{coro}\par
Formula (1) of Theorem \ref{primointbis} is very useful in applications,
expecially when $\LL=\LL^{**}$. \ In the next theorem we prove that  in the almost Gorenstein case the
blowing up  
  along a reflexive ideal $I$ is reflexive; this is not always true  (see Example
\ref{nonrifl}). Nevertheless, in  Example
\ref{anegativo}   we have $R$  almost Gorenstein, $\LL$
reflexive, but
$I$ not reflexive. \par First we recall the following property (see \cite{oz1}, Corollary 3.15).
\begin{(*)}\label{refl} \ \\
  \ Let $R$ be almost Gorenstein and let $J$ be a fractional ideal not
isomorphic to $R$, then  \vspace{-0.2cm}
$$\vspace{-0.2cm} J  \ \   is\
reflexive  \quad if\ and\ only\ if  \quad J:J\supseteq
R:\m.$$
\end{(*)}

 \begin{thm}  \label{corooz}   \ \\
Suppose 
$R$  almost Gorenstein and let  \ $\LL=\LL(I)$. Then
\begin{enumerate}
\item[(1)]  The equivalent conditions of the groups $({\rm A}), (\rm B)$ of Proposition
\ref{scopp.riflessivo} are equivalent to the following ones:
\begin{enumerate}
\item[ {\rm (C)}] \quad ${\rm (C_1)} \quad 
   \LL \supseteq R:\m.  \qquad {\rm (C_2)}\quad   I^{\nu}$ is reflexive. \\ \indent\quad ${\rm (C_3)}\quad\     I^n$ is
reflexive \  $\forall \ n \geq \nu$. \\ \indent\quad ${\rm (C_4)}\quad  \ I^n$ is reflexive
for some \ $  n \geq \nu$.
\end{enumerate}
\item[(2)]   If \ $\!   I\!\! \ $  is reflexive, then the equivalent
conditions $({\rm A}), (\rm B),({\rm C})$
hold, in particular \  $\LL$  is reflexive.
\end{enumerate}
\end{thm}
{\it Proof} \\
(1)  The equivalence of conditions (C)  is immediately achieved  by using
\ref{refl}.\par    
In order to prove the equivalence  (A) {\it if  and  only if} (C), we note
that  \par $\w I^n=(I^n)^{**}$
\ by Proposition \ref{Matlis2}, \ hence    \vspace{-0.2cm}
$$\vspace{-0.2cm} \w I^n= I^n
\quad if\ and\ only\ if  \quad I^n=(I^n)^{**}.$$
(2)  By applying as before Proposition \ref{Matlis2} we deduce that \ $I=
I^{**}=\w I$; \par
but this   is equivalent to \  \  $\w \subseteq I:I \subseteq \LL$.

\section{ Blowing up along the maximal ideal.}\vspace{-0.3cm}
 Our purpose   is  now to consider the  special  case  \ $I=\m$.  We denote by  $\LL\!:=
\!\LL(\m)$   {\it the  blowing-up} \  of $R $ along the maximal ideal,
  $e $   the   {\it  multiplicity}, \
 $\mu:= l_R(\m/\m^2) $    the
 {\it  embedding dimension}, $r$  the  {\it  Cohen-Macaulay type } of $R$; \ 
$x\in\m$ \ is    a  {\it  minimal reduction}   of  $\m$. \vspace{2mm}\\
\indent When $e=\mu$, namely $\m$ is {\it stable},  we can prove that the Gorensteiness
of the blowing up
$\LL$ is equivalent to the almost-Gorensteiness of the
ring $R$.  \par When  $e=\mu+1$, we get an explicit formula for 
  the length of the module $  R:\LL/\m^{\nu}  $. It turns out that this
length is zero if and only if $R$ is
Gorenstein and $\nu=2$. \par  In the cases
$\nu=2$ and
$\nu=3$ we state formulas involving the conductor $R:\LL$ \ which extend
some results of Ooishi valid for    Gorenstein
rings (see \cite{oo1}). \par
We begin with two simple remarks, useful in the sequel.    

\begin{rem} \label{primaug}   \
\begin{enumerate}
\item[(1)] \  $ l_R(R/R:\LL)=l_R(x(R:\m)/R:\LL)+(e-r)$.
\item[(2)] \ If $R$  is almost Gorenstein, then \  $  \LL=\LL^{**}$.
\end{enumerate}
\end{rem}
  {\it Proof}   \begin{enumerate}
 \item[(1)]  \ We know that $x\LL=\m\LL$ by property (1) of \ref{defLL}.
Therefore the inclusion \
$\m \subseteq x \LL$ \ implies that \ $ R:\LL \subseteq x(R:\m) \subseteq R
\subseteq  R:\m$. \ From this chain we get the
thesis.
\item[(2)] \ This is true by Theorem \ref{corooz}, since $I=\m$ is reflexive.
\end{enumerate}

\begin{(*)} \ { \bf CASE \    $e=\mu$.}\label{e=mu} \\
{\rm We recall that  \ $\!\m\! $ \ is said to be {\it stable} if \ 
$\LL=\m:\m$.  We have  the following
well known equivalent conditions for the {\it stability} of $\m$ (see
\cite{m}, Theorem 12.15):\par
  (1) \quad $\m$  {\it is stable} \par (2) \quad    $ e=\mu $ \par (3)
\quad     $\rho = e-1 $ \par (4) \quad     $
r=e-1   $.  } \end{(*)}

 \begin{prop} \label{casoe-1} \ \\
If \ $\m$ is stable, then    the following facts are
equivalent:\begin{enumerate}\item[(1)] \ R \ is almost Gorenstein
\item[(2)]
\ $\LL$ \ is Gorenstein.\end{enumerate}
 \end{prop}
{\it Proof}\\ By hypothesis \ $R:\LL=\m$ and $\nu=1$.  Hence if $R$ is
almost Gorenstein, then $\LL$ is Gorenstein by
Theorem
\ref{primointbis}.
 Vice versa, the hypothesis \ $ \LL=\m:\m$ \ implies that \ $c_{\LL}=c-e$.
\ Thus if $\LL$ is Gorenstein, then  condition
(2),(b) of Proposition
\ref{modcanLL} is satisfied and  $R$ is almost Gorenstein because\\
\centerline{ $2\delta
-c=2\rho-e=e-2=r-1 $.}

 \begin{(*)} \ \label{e=mu+1} { \bf CASE \
    $e=\mu+1$.} \\
 {\rm If $e=\mu+1$, the structure of $R$ is quite well understood, see e.g.
\cite{s}.  From the form of the {\it h-}polynomial    
\  $h(z)=1+(\mu-1)z+z^{\nu}$, \ one can infer that
  \ $\rho= \mu-1+\nu$. \  Moreover there are two possibilities depending on
the Cohen-Macaulay type $r$:\par
  (A):  \ If \  \ $\!  r <e-2$, \ then \ $\nu=2$;\par
  (B): \  If \ $r =e-2$, \ then \ $\m^2=x\m+(w^2)R,$ \ with \ $w \in \m
\setminus \m^2$
\ and \indent\qquad\quad \qquad$\m^3 \subset x \m$ \ (see \cite{s},
Prop. 5.1).} \end{(*)}

 We begin with a technical lemma.

\begin{lemma} \label{lemmae=mu+1} \quad
Assume that \  $r=e-2$. \  Then there exists an element $w\in
\m$ with \ $v(w)-e \notin v(\m:\m)$ \ such that:
\begin{enumerate}
\item[(1)] \ $\m\ =x(\m:\m)+w R $  \ and \ $w\m \subset x(\m:\m)$. \par
\item[(2)]  \ $\m^j=x \m^{j-1}+w^j R =$\par
\ \ \ \ \  \    $\! =x^{j-1} \m + x^{j-2} w^2 R+ .....+x w^{j-1} R + w^j R, \ \
\forall \ j=2,....,\nu$.
\item[(3)] \ $\m^3\subseteq x\m$.
\item[(4)] For every   element\
 $s \in \m:\m$ \ such that \ $v(s) >0$ \ we have
\vspace{-0.2cm} $$ \vspace{-0.2cm} sw^j\in x^{j-1}  \m, \ \ \forall \ j=2,....,\nu.$$
\end{enumerate}
\end{lemma}
{\it Proof} \begin{enumerate}
 \item[(1)]   The assumption \ $r=e-2$ \ means that \
$l_R(\m/x(\m:\m))=1$, \ hence by
\ref{jager}  there exists an element
$w\in
\m$ \ such that \ $v(w)-e \notin v(\m:\m)$ \  and \ $\m=x(\m:\m)+w R$. \ To
prove the inclusion \ $w\m \subseteq x(\m:\m)$ it suffices to consider the chain $x(\m:\m)\subseteq x(\m:\m)+w\m
\subset \m$. 
\item[(2)]We prove our claim by induction on $j$. Suppose \ $j=2$. \
From (1)
we have that  \ $ \m^2 \subseteq x \m +w\m = x \m
+w(x(\m:\m)+w R)\subseteq   x \m +w^2 R \subseteq \m^2$. \ Suppose now the
assertion true for $j$. Claim: \par
 $\m^{j+1}=x \m^{j }+w^{j+1} R  $\par
\ \ \ \ \ \ \ \  $=x^{ j } \m + x^{j-1} w^2 R+ .....+x w^{j } R + w^{j+1}
R $.\\
By using repeatedly the inductive hypothesis we get \\
 $\m^{j+1}=x \m^{j }+w^j\m \subseteq x \m^{j }+w \m^j =  x \m^{j }+w (x
\m^{j-1}+w^j R )$  \par  $ \ \ \ \ \ \ \ \ \subseteq x \m^{j }+w^{j+1} R
\subseteq \m^{j+1}$.\\
We are left to prove the second equality of the claim. We have:\\
\ $\m^{j+1}=  x^{j-1} \m^2 + x^{j-2} w^2 \m+ .....+x
w^{j-1}
\m + w^j
\m$ \par
$ \ \ \ \ \ \ \ \ \subseteq  x^{j-1}
\m^2 + x^{j-2} w \m^2+ .....+x w^{j-2} \m^2 + w^{j-1} \m^2 $ \par
$ \ \ \ \ \ \ \ \ = x^{j-1} (x \m
+w^2 R ) +
x^{j-2} w  (x \m +w^2 R ) + .....+x w^{j-2}  (x \m +w^2 R ) $ \par $ \ \ \ \ \ \ \ \ \ \ \ \ + \ w^{j-1}  (x
\m +w^2 R ) $ \par $ \ \ \ \ \ \ \ \ \subseteq  x^j \m +
x^{j-1} w^2 R+ .....+ x^2 w^{j-1} R+x w^{j } R + w^{j+1} R \subseteq
\m^{j+1}$. \\ For the last but one inclusion we
have used the fact that \par
$ \ \ \ x^{j-1} w   \m  + .....+x^2 w^{j-2}  \m  +  x w^{j-1} \m$ \par $ = x^{j-1} w
(x(\m:\m)+w R)  + .....+    x w^{j-1} (x(\m:\m)+w R) $ \par
$\subseteq x^{ j } \m +
x^{j-1} w^2 R+ .....+x^2 w^{j-1} R +x w^j R$.
\item[(3)]  As seen in the proof of item (2), $ \m^2 = x \m +w\m$. \ Hence  \\ $ \m^3 = x \m^2 +w\m^2 \subseteq x\m,
$ because \ $w\m^2\subseteq x\m$ by item (1).
 \item[(4)]   
Let \ $s\in\m:\m$\ be such that   \ $v(s) >0$. \ We proceed by induction on $j$. Suppose   $j=2$. \ By
item (1)   there exist  \ $y\in \m:\m $ and $ \  a\in R$
\ such that \
$sw=xy+a w$. \ If  \ $v(a)=0$, \ then \ $v(s-a)=0$, \
contradicting the fact that \ $v(w)-e \notin v(\m:\m)$. \ Hence \ $a\in\m$.
\ Thus \   $sw^2 =xyw+aw^2\in
x\m$, \ because \ $\m^3 \subseteq x\m$.
\ Assume now the inductive hypothesis
\ $\displaystyle\frac{sw^j}{x^{j-1}} \in m$, \ then \
$\displaystyle\frac{sw^j}{x^{j-1}}=xz+b w$, \ with \ $z\in \m:\m,
\  b\in R$, \ i.e., \ $(\displaystyle\frac{sw^{j-1}}{x^{j-1}}-b)w= xz$. \
Since  the element \ $
 \displaystyle\frac{sw^{j-1}}{x^{j-1}}  $ \ has a positive valuation, by
the same reasoning  as above we conclude  that \
$b\in \m$.
\ Therefore \ $ \displaystyle\frac{sw^{j+1}}{x^{j-1}}= xzw+bw^2\in x\m$, \
which is our thesis. \end{enumerate}

\begin{prop} \label{prope=mu+1}  \    \begin{enumerate}
\item[(1)]   If \ \ $r=e-2$, \ then \quad  $ e=\mu+1$.
\item[(2)] If \ \ $e=\mu+1$, \ then\quad
 $l_R(x(R:\m)/R:\LL)=1.$ \end{enumerate}
\end{prop}
{\it Proof}
\begin{enumerate} \item[(1)] By Lemma \ref{lemmae=mu+1},(2), \ $\m^2=x\m
+w^2 R$. \ Hence
\vspace{-0.2cm} $$\vspace{-0.2cm}l_R(\m/\m^2)=l_R(\m/x\m)-l_R(\m^2/x\m) =
e-1.$$
\item[(2)]  We shall prove
that the
$R$-module
\vspace{-0.2cm} $$\vspace{-0.2cm}   R:\m/x^{\nu-1}(R:\m^{\nu}) \simeq
x(R:\m)/R:\LL$$
\ is monogenous generated by $\overline1$. We divide the proof  in two
parts, following cases \ (A), $r <e-2$ \ and \ (B),  \ $r =e-2$ above. \\
Case (A).  \  $\nu=2$. \   Since \ $\m^2= x\m+(a)R, \ a \notin x\m$, \ we have
that \ $x(R:\m^2)= (R:\m)
\cap \Big ( \displaystyle  \frac{x}{a}\Big )R$. \ If \ $y\in R:\m$, \  then
\ $ya
\in \m^2$ \  and we can write \ $ ya= xr+as,$ \ with $ \ r\in \m, \ s\in
R$, \ namely $ \  y=
\displaystyle  \frac{x}{a} r+s$, \ so $ \ \overline y = \overline 1 s$.\\
Case (B). \  We want to prove that if \ $s \in \m:\m$ has a positive
valuation, then \
$s \in x^{\nu-1} (R:\m^{\nu})$. \ By item  (2) and (4)  of
Lemma
\ref{lemmae=mu+1} we have \vspace{-0.2cm} $$\vspace{-0.2cm}
\displaystyle\frac{s\m^{\nu}}{x^{ \nu-1 }} \in s\m + s
\displaystyle\frac{w^2}{x} R +.... + s
\displaystyle\frac{w^{\nu-1}}{x^{ \nu-2 }} R + s
\displaystyle\frac{w^{\nu}}{x^{ \nu-1 }} R \subseteq \m.$$\end{enumerate}

\begin{thm} \label{teo,e=mu+1}    \ \\  Let $e=\mu+1$.\  Then
$$l_R(R:\LL/\m^{\nu}) =r-1+(e-1)(\nu-2).$$
\end{thm}
 {\it Proof }\\ We have to compute the difference \ $l_R(R/\m^{\nu}) -
l_R(R/R:\LL)$. \ As recalled in \ref{e=mu+1} \
$e\nu-l_R(R/\m^{\nu})=\rho= e-2+\nu$. \  Combining the above results
\ref{primaug} and
  \ref{prope=mu+1} we obtain \ $l_R(R/R:\LL)= e-r+1$.  \ The conclusion
follows.

\begin{coro} \label{coroe=mu+1}    \ \\  Let $e=\mu+1$. \ Then
\begin{enumerate}
\item[(1)] \ $ R:\LL=\m^{\nu}$ \quad if and only if \quad  $R$  is   Gorenstein
and \ $\nu=2$.
\item[(2)] \  $ \sum_{i\notin \Gamma,i \neq 1} \ (r_i-1)=
d(R:\Lambda)+l_R(\LL^{**}/\LL)+(\nu-2)$.
\item[(3)] \  $R$ is almost Gorenstein  \quad if and only if \quad   $
\nu=2 $ \ and
\ $\w\LL=\LL$.
\end{enumerate}
\end{coro}
 {\it Proof }
 \begin{enumerate}
\item[(1)] It follows directly from Theorem \ref{teo,e=mu+1}.
\item[(2)] As recalled in \ref{e=mu+1} \   $\rho= \mu-1+\nu$, \
then   the formula of Theorem \ref{primoint} 
gives:   \vspace{-0.2cm}
$$\vspace{-0.3cm}l_R(R:\LL/\m^{\nu}) =(\mu+1)\nu -2(\mu-1+\nu)+\sum_{i\notin \Gamma} \
(r_i-1)- d(R:\Lambda)-l_R(\LL^{**}/\LL).$$
By comparing with Theorem \ref{teo,e=mu+1} the thesis follows.
\item[(3)] This is clear after observing that  item (2) combined with 
equality (1)
of  Proposition \ref{dopodefc}  becomes:
\vspace{-0.2cm}
$$\vspace{-0.2cm}(2\delta-c)-(r-1)=l_R(\w\LL/\LL)+(\nu-2).$$\end{enumerate}
\begin{coro} \ \\  
    \  $r=e-2$  \ \ and \ $ R:\LL=\m^{\nu}\quad  if\ and\ only\ if\quad R\ \  is \  Gorenstein\
with\ \ e=3.$ \end{coro}

 \begin{(*)} \ { \bf CASE \ $\nu=2$.} \label{nu=2}  \\
  {\rm We recall that in this case the invariants  $\rho,
e, \mu $ are related by the equality:\vspace{-0.2cm}
  $$\rho= 2e-\mu-1.$$ }
\end{(*)}

\begin{prop} \label{s=2}  \ \\ Assume \ $ \nu=2$.
\begin{enumerate}
\item[(1)] \ $ 2e+r l_R(R:\Lambda/\m^2) \leq ( r+1)(\mu+1).$
\item[(2)]  If $R$ is almost Gorenstein, then \ $
 e-(\mu+1)= \displaystyle\frac{r-1}{2}- \displaystyle
\frac{1}{2}l_R(R:\LL/\m^2)$.\par
 In particular: \item[]  if  $R$ is  Gorenstein, then   $  e=\mu+1$ \ and \
$ R:\Lambda=\m^2;  $
\item[] if  $R$  is a  Kunz ring (namely almost Gorenstein of   type
2),   then  \vspace{-0.2cm}  $$\vspace{-0.2cm}  e=\mu+1  \ and \
  l_R(R:\LL/\m^2) = 1.$$    
\end{enumerate}
 \end{prop}
 {\it Proof }
 \begin{enumerate}
\item[(1)]   The  inequality
 follows directly from Corollary
\ref{minor}.
\item[(2)] This is Theorem \ref{primointbis} with \ $ \nu=2$,  $\rho=
2e-\mu-1$ and
$\LL=\LL^{**}$.
\end{enumerate}

\noindent We deduce from   Proposition \ref{s=2} that  if
$R$ is almost Gorenstein, then:
\vspace{-0.2cm}
$$(*) \ \ \ \ \ \   \ R:\LL=\m^2 \ \ \  if\  and\   only\   if\ \ \
e-(\mu+1)=
\displaystyle\frac{r-1}{2} \ \ {\it and} \ \ \nu=2. \vspace{-0.2cm} $$
This equivalence was already known for Gorenstein rings:  assertion  $(*)$
in the
case $r=1$ is exactly Corollary 7 of \cite{oo1}. \\
We remark that  there
exist almost Gorenstein rings satisfying the condition
\centerline{$e-(\mu+1)=
\displaystyle\frac{r-1}{2}   $}\\
 \ with \ $\nu \neq2$ \ (see Example \ref{ninon2}).
The next corollary shows that this  cannot happen when  $R$ is Gorenstein.

\begin{coro} \label{caseGor}      \ \\ Let  $R$ be a   Gorenstein ring.
Then the following conditions are equivalent:\begin{enumerate}
\item[(1)]
  \quad $ e= \mu+1.$ \item[(2)]\quad $ \nu=2. $ \item[(3)]\quad $
R:\LL=\m^2. $ \end{enumerate}
 \end{coro}
{\it Proof} \\
 If $ \ e=
\mu+1,$ we have  that \ $\nu=2$ \ by \ Corollary  \ref{coroe=mu+1}. To
conclude the proof it   suffices to apply
Proposition
\ref{s=2}.   \\[0.25cm] 
\indent In the case    $R$   Gorenstein, Corollary \ref{coroe=mu+1} combined with
Corollary \ref{caseGor} gives  Proposition
12 of  \cite{oo1}:

\centerline{$\gamma_R=\m^2 $   \ {\it if and only if \ \  $ \LL=\R$ \ and} \
$  e=\mu+1$.}
\noindent The following proposition states a more general relation between
$  e $ and   $\mu+1$.

\begin{prop} \label{dopooishi} \ \\
 The following conditions are equivalent:\begin{enumerate}
\item[(1)] \ $\gamma_R=\m^2$
\item[(2)] \ $ \LL=\R, \   e-(\mu+1)= \displaystyle
\frac{1}{2}(2\delta-c)$ \  and \  $\nu=2$.\end{enumerate}
\end{prop}
\vspace{-0.2cm} {\it Proof}\\
 If \  $ \gamma_R=\m^2, $ \ then we have   \ \ $ \m^2:\m^2=\R=\LL $ \ and
$ \ \nu=2; $  \quad hence we get \\
 \centerline{$c=2e$ \quad and \quad $2 \rho= 2e+2\delta-c$,} because
condition (c) of
Proposition \ref{modcanLL} is verified.  On the other hand, since \
$\nu=2$, \ $\rho=2e-\mu-1$. \ By comparing the two
equalities we see that (2) holds.\\
Conversely, \ $\gamma_R= R:\LL$ \ and \  $\delta= \rho=2e-\mu-1 $ \ imply
that \\ \centerline{$ 2 \delta= 4
e-2\mu-2= 2e-2(\mu+1) +c  $.} \
Hence $ \  c=2e$, \  and   again by Proposition
\ref{modcanLL} \ we obtain that
$\gamma_R=R:\LL=\m^2$.

 \begin{(*)} \ { \bf CASE \   $\nu=3$.}\\
 {\rm This case has been considered by Ooishi in \cite{oo1}. Statement  (3)
of the next proposition extends to   almost Gorenstein rings Proposition 8
of his quoted paper, valid  in the case
$r=1$. } \end{(*)}

\begin{prop} \label{props=3}   \ \\
 Assume that $\nu=3$.
      \begin{enumerate}
\item[(1)] If \ $r=2,$ \ then \ $  e-(\mu+1)+\displaystyle\frac{2}{3}
l_R(R:\Lambda/\m^2)  \leq
l_R(\m^2/\m^3)$
\item[(2)] If the h-polynomial is symmetric,  \ then \vspace{-0.2cm} $$
l_R(R:\LL/\m^3)
\leq  3
\mu
\displaystyle \frac{(r-1) }{r}. \vspace{-0.2cm}$$
\item[(3)] If $R$ is almost Gorenstein, then
  the   h-polynomial is  symmetric   if   and  only  if
  \ $ l_R(R:\LL/\m^3) =r-1 $ \ and \ $ e = 2 \mu.$
\end{enumerate}
 \end{prop}
 {\it Proof}  \begin{enumerate}
\item[(1)] It suffices to apply Corollary \ref{minor}.
\item[(2)] Of course,   $h_{\m}(z)=
1+(\mu-1)z+a_2z^2+a_3 z^3$  is symmetric if and only if  $  a_3=1$ \ and \
$   a_2=\mu-1$. \ In this case we obtain
\
$ e=  2
\mu$ \  and  \
$\rho = 3 \mu$.
\ The inequality follows from Corollary \ref{minor}.
\item[(3)]This comes directly from the main formula of Theorem
\ref{primointbis}.  \end{enumerate}

 Using
again the formula of Theorem \ref{primointbis} we get immediately the next result.

\begin{coro} \label{ooishi}    \ \\
 Suppose $R$ almost Gorenstein. If
\  $e  = 2 \ \mu $, \ then \\
 \centerline{$R:\LL=\m^{\nu}  \quad $ if   and  only  if\ \ \quad  $\rho=
\nu \  \mu + \displaystyle\frac{r-1}{2}.$}
\end{coro}

We notice that in the case \ $r=1$ and $\ \nu=3$ \  result   \ref{ooishi} gives
again Proposition 8   of \cite{oo1}; however
  there exist almost Gorenstein, not Gorenstein, rings satisfying the 
equivalent conditions of \ref{ooishi}   (see Example
\ref{eduemi}).

 \section{Examples.}\vspace{-0.3cm}

In all examples listed below we suppose that \ $R=\cx[[t^h]], \ h \in
v(R)$, is a semigroup  
ring and \ that $\LL=\LL(I)$ \ is the blowing-up of $R$ along the specified ideal $I$.\vspace{3mm}  \\
\noindent Notation $< ... >$ means "the semigroup generated by ... ".
Notation $a-b$ in
the semigroup means "all the integers between $a$
and $b$". Notation  $a\rightarrow$   means "all the integers   $\geq a$".

\begin{example}\label{negativo} {\rm  (See the first remark in
\ref{nonaltro}).\\ Let \
$v(R)=\{0,10,12,20  \rightarrow\}$.
\begin{enumerate}
\item[(1)] If $\  I=(t^{10},t^{12}) $, \ then  \ $v(\LL)=<2,21>$, \ hence  \
$c=c_{\Lambda}= 20 $ \   and \
$ l_R(R/\gamma_{R})-l_R(\LL/\gamma_{\LL}) = -\rho $.
\item[(2)] If
$\   I=\m $, \ then \ $v(\LL)=<2,11> $, \ hence \   $c-c_{\Lambda}=10=e$ \
and \ $ l_R(R/\gamma_{R})-l_R(\LL/\gamma_{\LL})= e-
\rho=-2$.
\end{enumerate}}
\end{example}

\begin{example}\label{cinque} {\rm(See the second remark in
\ref{nonaltro}).\\  Let \
$v(R)=  \{0,5,10,11,12,15,16,17, 19\rightarrow\}$, i.e., $v(\m)= <5,11,12,19>$,  \ and let \  $ I=\m $.  
 \ We have \ $\nu=2$ \ and
\ $v(\Lambda)=\{0,5,6,7,10 \rightarrow\}$. \  Moreover:  \\ $ v(R:\m^2)=\{0,5,6,7,9
\rightarrow\} $, \ hence \ $R:\m^2 \subset\R$, \
but \ $ c-c_{\Lambda}=9  < e \nu=10$.} \end{example}

\begin{example} \label{ninon2}
 {\rm  (See the third remark in \ref{nonaltro} and also the remark after
\ref{s=2}).\\ Let \  $v(R) =
\{0,7,8,12,13,14,15,16,18
\rightarrow\}$, i.e., $v(\m)=<7,8,12,13,18> $.
\\
 $R$ is almost Gorenstein and its  {\it h-}polynomial  is \ $ h(z)=1+4 z+z^2+z^4$.
\\ We
have \  $\m^4= t^{28}\R $, \ hence \ $\LL= \R$ \
and \ $\m^4 \subset  R:\LL = \gamma_R= t^{18} \R$. \\  In this case formula
\  $ \gamma_R =(R:\Lambda)
\gamma_{\Lambda} $  \ holds, but  \ $c-c_{\LL}=18 < e \nu= 28$.\\  Since \
$e=7, \mu=5, r=3$,
\   condition \
$e-(\mu+1)=
\displaystyle\frac{r-1}{2}$ \  is satisfied, but $\nu=4$. }
\end{example}

\begin{example} \label{anegativo}   {\rm (See the fourth  remark in
\ref{nonaltro} and also the remark after
\ref{sottocaso}).\\Let \
 $v(R) = \{0, 5, 10, 15, 20, 21, 25, 26, 30-32, 35-37,40-
42, 45- 48, 50- 53, 55-
58, 60 \rightarrow\} $, i.e., $v(\m)=<5,21,32,48>$. \\  $R$ is an almost Gorenstein ring with
Cohen-Macaulay type
$r= 3$  and
$ e= \mu+1$.
\begin{enumerate}
\item[(1)] If $\  I=\m $, \ then  \ $v(\LL)= \{0, 5,10, 15,16, 20, 21, 25-27,
30- 32, 35- 37,
40-43, 45- 48, 50\rightarrow\}$ \ and $\nu=2$, \ hence  \ $c-c_{\Lambda}=
10=e\nu $ \ and $\LL$ is reflexive, \  but \  $\m^2
\subset  R:\LL $.
\item[(2)] If \
$  I= (t^{31}, t^{32}, t^{40}) $, \ then \ $I$ is not reflexive, whereas \
$\LL =I^4:I^4= \R $ \   is reflexive.
\end{enumerate}}
 \end{example}

\begin{example} \label{cnonnullo}   {\rm (See   the remark after
\ref{primoint}).\\
 Let $   v(R) =\!\{0,\!10,\!20,21\!,25\!,26\!,30-36,40-47,50\rightarrow\!\}$, \\i.e., $   v(\m)= 
<10, 21, 25, 26,32,33,34>$,\ and let \
$I=\m$.\\
Here \    
$v(\LL)=\{0,10,11,15,16,20-27,30\!\rightarrow\}$,  
  $\LL\!=\! \LL^{**}\! \!=\! \m^2\!:\!\m^2$    and\   $  R:\LL\! =\m^2$. \\
Since the type sequence of $R$ is \
$[3,2,1,2,1,3,1,1,1,1,1,1,2,1,1,1,1,1,1,1,2], $ \  \ \ \  $v(R:\Lambda)=
\{20,30,31,35,36,40-47, 50\rightarrow\}, \ \Gamma=\{ 3,7,8,12\rightarrow\}$, \ we have that \vspace{-0.2cm}
$$\vspace{-0.2cm}\sum_{i \in \Gamma}r_i=  15  \quad {\rm and }\quad    d(R:\LL)=
l_R(\R/\LL)-\displaystyle\sum_{i
\in
\Gamma}r_i=17-15= 2.$$ }
 \end{example}
\begin{example} \label{GnonGor} {\rm(See   the remark after
\ref{minorbis}).\\ Let \ $ v(R) =\{0 ,11,12,15,22-
27, 29, 30, 33- 42, 44\rightarrow\}$, \\ i.e., $ v(\m)=<11,12,15,25,29>$, \ and let \
$I=\m$. \\
$R$ is a Gorenstein ring and its {\it h-}polynomial \  $h_I(z)=1 + 4z  + 2z^2 +
2z^3 +
2z^4$ \ is not symmetric.
 \ We have $\nu=4$ and $\rho=22$, \ hence \ $2 \rho=e \nu$. }
 \end{example}

\begin{example} \label{comeGnonGor}   {\rm (See   the remark after
\ref{minorbis}).\\Let $R$ be such that \
 $v(\m):= <10,23,55,58,82>$ \ and let \ $I=\m$.\\
$R$ is almost Gorenstein with Cohen-Macaulay type $r=3$ and its
{\it h-}polynomial \  $h_I(z)=
1+4z+z^2+2z^3+2z^4$ \ is not symmetric. \ We have  $
\rho=20, \
\nu=4$. \ Hence $R$ verifies the condition \ $2 \rho= e
\nu$.}
\end{example}

\begin{example} \label{condverif} \ \\  {\rm In this example $R$ is
an almost Gorenstein ring with Cohen-Macaulay type $r=3$, verifying
the equivalent conditions of Theorem \ref{primointbis}.\\
Let \
$v(\m):=<10,16,95,99>$  \
and let \ $I=\m$.\\
Its {\it h-}polynomial is   $h_I(z)=1+3z+2z^2+2z^3+2z^4$ and \ \ $c=124$.    Since  $\rho=21$
and $\nu=4$, \
 we have \ $ 2\rho= e \nu+r-1$. \ It follows that  $\LL= \LL^{**}$ is
Gorenstein and $c_{\Lambda}= c- e
\nu= 84$.}
\end{example}

\begin{example} \label{nonrifl}  {\rm (See the remark after
\ref{sottocaso}).\\Let $R$ be such that \
$v(\m):= < 6,11,16,20,25>$ \ and let \  $I = \m$.\\
The blowing up \ $ \LL =  \m^2:\m^2$ \ is not reflexive. }\end{example}

\begin{example} \label{eduemi}  {\rm (See the remark after
\ref{ooishi}).\\
 Let  \ $v(R) = \{0,8,10,13,15,16,18,20,21,23-26,28\rightarrow\}$, \\i.e., $v(\m)=<8,10,13,15>$, \ and let
\  $I = \m$.\\
$R$ is an almost Gorenstein ring  with Cohen-Macaulay type $r=3$, verifying
the conditions of
Corollary \ref{ooishi}. In fact,  
  $e= 2 \mu$  and the {\it h-}polynomial is
\
$h_I(z)= 1+3 z+2z^2+2z^3$, \ hence
\ $\rho= 13= \nu \mu+(r-1)/2$. \ Thus \ $R:\LL=\m^3$ \ and by
\ref{primointbis} \ $c_{\LL}=c-e\nu=4$.}\end{example}

{\small {\it  A. Oneto }  \\ Dipem,  
  Universit\`a di Genova,   P.le Kennedy, Pad. D -  I 16129 Genova (Italy); \\oneto@dipem.unige.it\\
  { \it E. Zatini} \\ Dima,  Universit\`a di Genova,  Via Dodecaneso 35 -  16146 Genova (Italy);\\
zatini@dima.unige.it}

\end{document}